\documentclass[12pt]{article}
\begin{document}

\title{Complementary Polynomials From Rodrigues' 
Representations For Confluent And Hypergeometric 
Functions And More}  
\author{H. J. Weber\\Department of Physics\\
University of Virginia\\Charlottesville, 
VA 22904, U.S.A.}
\maketitle
\begin{abstract}
Complementary polynomials of Legendre 
polynomials are briefly presented, as well 
as those for the confluent and hypergeometric 
functions, relativistic Hermite polynomials 
and corresponding new pre-Laguerre polynomials. 
The generating functions are all given in closed 
form and are much simpler than the standard ones. 
Some are simply polynomials in two variables.   
New recursions and addition formulas are derived.
\end{abstract}
\leftline{MSC: 33C45, 33C55, 34B24}
\leftline{Keywords: Rodrigues' formula, complementary  
polynomials}  

\section{Introduction}

The novel method introduced in~\cite{hjw} 
transforms the Rodrigues' formula of given 
polynomials step by step thus defining the 
complementary polynomials recursively. The 
most useful property of complementary 
polynomials is that their generating 
function can always be given in closed form. 
As a rule, it is simpler than the standard 
generating function of the polynomials 
defined by the Rodrigues' formula.  

Hence polynomials satisfying a Rodrigues' formula 
are accompanied by their complementary polynomials. 
The classical polynomials that are important in 
mathematical physics are such cases: Legendre 
polynomials and the polynomial components of 
associated Legendre functions form such pairs; 
so do Laguerre and associated Laguerre polynomials, 
etc. Some are discussed below, along with with new 
and old results (with new simple derivations).      

For associated Laguerre polynomials, 
it leads to so simple a generating function in 
closed form that the H-atom's wave functions 
become truly accessible to undergraduates. 
This is now implemented in Chapt.~18 on 
pp.~892-894 of the textbook~\cite{awh}. 

For complementary Legendre polynomials, which 
are the polynomial parts of associated Legendre 
functions~\cite{trio}, the generating 
function is simply a polynomial in two variables. 

For Romanovsky polynomials their complementary 
polynomials are again Romanovski polynomials. 
This furnishes a generating function for them 
in closed form~\cite{rom}. These polynomials 
have recently played a prominent role in 
physics~\cite{rwak}.   

Here we also consider generalized Rodrigues' 
formulas without the frame-work of a 
hypergeometric-type ODE and Pearson's ODE 
to make its operator self-adjoint.  
Yet, they generate polynomials along 
with their complementary ones that 
have a closed-form generating function.   

The normalizations of complementary 
Legendre polynomials are revisited in 
Sect.~2. The generating function and some  
applications of the confluent hypergeometric 
ODE are dealt with in Sect.~3. The 
hypergeometric ODE and its complementary 
polynomials are the subject of Sect.~4. 
Relativistic harmonic oscillator polynomials 
are treated in Sect.~5. Similar pre-Laguerre  
polynomials are introduced in Sect.~6.  

\section{Complementary Polynomials For\\
Legendre Polynomials}

{\bf Definition~2.1} Complementary 
polynomials ${\cal P}_\nu(x,l)$ are defined 
recursively by Rodrigues' formula~\cite{trio} 
\begin{eqnarray*}
(-1)^l 2^l l! P_l(x)=\frac{d^{l-\nu}}{dx^{l-\nu}}
\left((1-x^2)^{l-\nu}{\cal P}_\nu(x,l)\right),~
\nu=0,1,2,\dots,l , 
\end{eqnarray*}
where $P_l(x)$ are the Legendre polynomials. 
 
{\bf Theorem~2.2} {\it The generating function 
has the closed form}  
\begin{eqnarray}
{\cal P}(y,x,l)=(1-y(1+x))^l(1+y(1-x))^l.  
\label{gef}
\end{eqnarray}    
{\it Proof.} This is Eq.~(12) in Sect.~3  
of~\cite{trio}.$~\bullet$ 

Thus, ${\cal P}(y,x,l)$ is simply a product 
of two forms each linear in $x,~y$ and raised 
to the lth power. In terms of the angular 
variable $x=\cos\theta,$ common in physics 
applications,   
\begin{eqnarray}
{\cal P}(y,\cos\theta,l)=(1-2y\cos^2\frac{
\theta}{2})^l(1-2y\sin^2\frac{\theta}{2})^l.  
\label{trig}
\end{eqnarray}
Binomially expanding Eq.~(\ref{trig}) yields 
Eq.~(14) in~\cite{trio} in trigonometric form.  

{\bf Corollary~2.3} {\it The basic composition 
laws} 
\begin{eqnarray*}
[(1-xy)^2-y^2]{\cal P}(y,x,l)={\cal P}(y,x,l+1), 
\end{eqnarray*}
\begin{eqnarray*}
{\cal P}_\nu(x,l+1)={\cal P}_\nu(x,l)-2x\nu{\cal 
P}_{\nu-1}(x,l)-(1-x^2)\nu(\nu-1){\cal P}_{\nu-2}(x,l) 
\end{eqnarray*}
{\it hold.} 

{\it Proof.} The generating function identity 
follows at once from Eq.~(\ref{gef}). When 
expressed in terms of complementary polynomials 
this gives the recursion relation.~$\bullet$   

For many other composition laws we refer to~\cite{trio}.  

Partial derivatives of the generating function, such as 
\begin{eqnarray*}
\frac{\partial{\cal P}(y,x,l)}{\partial y}
=-2l[x+y(1-x^2)]{\cal P}(y,x,l-1),
\end{eqnarray*} 
and 
\begin{eqnarray*}
\frac{\partial{\cal P}(y,x,l)}{\partial x}
=-2ly(1-xy){\cal P}(y,x,l-1),
\end{eqnarray*} 
yield the identity 
\begin{eqnarray}
y(1-xy)\frac{\partial{\cal P}(y,x,l)}{\partial y}
=[x+y(1-x^2)]\frac{\partial{\cal P}(y,x,l)}
{\partial x}. 
\label{legi}
\end{eqnarray}
  
{\bf Corollary~2.4} {\it The following 
three-term recursion is valid}
\begin{eqnarray*}
{\cal P}_\nu-2(\nu-1-l)x{\cal P}_{\nu-1}=
(\nu-1)(\nu-2l-2)(1-x^2){\cal P}_{\nu-2}. 
\end{eqnarray*}
{\it Proof.} Equation~(\ref{legi}), when 
expressed in terms of complementary polynomials, 
gives the recursion 
\begin{eqnarray*}   
\nu{\cal P}_\nu-\nu(\nu-1)x{\cal P}_{\nu-1}=
x{\cal P'}_\nu+\nu(1-x^2){\cal P'}_{\nu-1}. 
\end{eqnarray*}
Using the basic differential recursion~\cite{hjw} 
Cor.~4.2
\begin{eqnarray}
{\cal P}'_\nu(x,l)=-\nu(2l-\nu+1){\cal P}_{\nu-1}(x,l)
\label{bdrec}
\end{eqnarray}
allows us to eliminate the derivatives to obtain the 
3-term recursion.~$\bullet$ 

From a comparison of the ODE satisfied by the 
complementary polynomials with the ODE of the 
polynomial components ${\cal P}_l^m(x)$ of 
associated Legendre functions, one finds~\cite{trio} 
\begin{eqnarray*}
{\cal P}_{l-m}(x,l)=N_l^m {\cal P}_l^m(x);~m=0,1,\ldots,l.
\end{eqnarray*}

For even $l-m$ the normalizations $N_l^m$ (see Eq.~(43) of 
\cite{trio} and Eq.~(\ref{2.3}) below) follow from the 
special values 
\begin{eqnarray*}
{\cal P}_{2\nu}(0,l)=(-1)^\nu(2\nu)!\left(l\atop \nu\right)
\end{eqnarray*}
implied by the generating function ${\cal P}(y,0,l)=(1-y^2)^l$ 
at $x=0.$  For odd $l-m$ we now derive a new 
formula for the normalizations.  

Comparing the differential recursion~\cite{hjw} (Cor.~4.2) 
(see Eq.~(\ref{bdrec}) above) with the well known 
\begin{eqnarray}
{\cal P}_l^m=(-1)^m\frac{d^mP_l(x)}{dx^m}
=-\frac{d}{dx}{\cal P}_l^{m-1}(x),~m=0,1,\ldots,l
\label{2.1}
\end{eqnarray}
yields, when expressed in associated Legendre 
polynomials, 
\begin{eqnarray}
N_l^m{\cal P}_l^{'m}(x)=-(l-m)(l+m+1)N_l^{m+1}
{\cal P}_l^{m+1}(x). 
\label{2.2}
\end{eqnarray}
Comparing Eqs.~(\ref{2.1}) and (\ref{2.2}), we 
obtain Eq.~(\ref{2.4}) for even $l-m:$ 

{\bf Theorem~2.5} {\it The normalizations of 
complementary Legendre polynomials are given by}
\begin{eqnarray}
N_l^m=(-1)^m\left(l\atop (l-m)/2\right)
\frac{(l-m)!!(l-m)!}{(l+m-1)!!},~l-m~~\rm{even}
\label{2.3}
\end{eqnarray}
\begin{eqnarray}
N_l^{m+1}=\frac{N_l^m}{(l-m)(l+m+1)},~-l<m<l.
\label{2.4}
\end{eqnarray}
So, for $N_l^{m+1}$ the index $l-m-1$ is odd, and 
this simple recursion renders unnecessary the 
lengthy developments in \cite{trio}~Sect.~4 
after Example~4.1. 

\section{Confluent Hypergeometric ODE} 

We write the standard regular solution 
of the confluent hypergeometric ODE  
\begin{eqnarray}
xy''+(c-x)y'-ay=0
\label{confl} 
\end{eqnarray}
as $M(a,c,x),$ as usual. Then, in the 
notations of \cite{hjw}, 
\begin{eqnarray*}
\sigma(x)=x,~\tau(x)=c-x,~\sigma'=1,~\tau-\sigma'=c-1-x. 
\end{eqnarray*}
Thus, Pearson's ODE for the weight function $w(x)$ 
\begin{eqnarray*}
xw'=(c-1-x)w(x)
\end{eqnarray*}
is solved by $w(x)=x^{c-1}e^{-x},x\in [0,\infty).$ 

{\bf Definition~3.1} {\it The Rodrigues' formula 
for the polynomial solutions $R_l(x)$ of Eq.~(\ref{confl}) 
is~\cite{hjw} (Eq.~(3))}  
\begin{eqnarray*}
R_l(x)=x^{1-c}e^x\frac{d^l}{dx^l}(x^{l+c-1}e^{-x}),~
l=0,1,\ldots 
\end{eqnarray*}
{\it and for the complementary 
polynomials~\cite{hjw}~(Eq.~(5)) it is} 
\begin{eqnarray}
{\cal P}_\nu(x,l)=x^{1-c+\nu-l}e^x\frac{d^\nu}{dx^\nu}
(x^{c-1+l}e^{-x}),
\label{comcon} 
\end{eqnarray} 
{\it so that $R_l(x)={\cal P}_l(x,l).$ The $R_l$ 
also depend on the parameters of the ODE, but 
these are suppressed for simplicity.} 

{\bf Lemma~3.2} {\it The ODE for complementary 
polynomials is} 
\begin{eqnarray}
x\frac{d^2{\cal P}_\nu}{dx^2}+(l-\nu+c-x)
\frac{d{\cal P}_\nu}{dx}=-\nu{\cal P}_\nu(x,l), 
\label{3.2}
\end{eqnarray}
{\it so that}
\begin{eqnarray*}
&&{\cal P}_\nu(x,l)={\cal P}_\nu(0,l)M(-\nu,l-\nu+c,x),
\\&&R_l(x)={\cal P}_l(0,l)M(-l,c,x),~{\cal P}_\nu(0,l)
=\nu!\left(c-1+l\atop \nu\right). 
\end{eqnarray*}
{\it Proof.} The ODE~(\ref{3.2}) follows from~\cite{hjw} 
Theor.~5.1. Since it is the confluent hypergeometric ODE 
again, this implies the other relations.~$\bullet$

Thus, all polynomial solutions of the confluent 
hypergeometric ODE are summed in the following 
1-parameter generating function in closed form, 
which is quite simple given the complexity of the 
regular solution of the confluent hypergeometric 
ODE.  

{\bf Theorem~3.3} {\it The generating function 
of the complementary polynomials takes the form}  
\begin{eqnarray}
{\cal P}(y,x,l)=\frac{[x(1+y)]^{c-1}e^{-x(1+y)}}
{x^{c-1}e^{-x}}(1+y)^l=(1+y)^{c-1+l}e^{-xy}. 
\label{confgenf}
\end{eqnarray} 
\begin{eqnarray*}
{\cal P}_\nu(x,l)=\sum_{\mu=0}^\nu\left(c-1+l\atop 
\nu-\mu\right)\frac{\nu!}{\mu!}(-x)^\mu.
\end{eqnarray*}
{\it Proof.} The generating function follows from 
\cite{hjw} Theor.~3.2 by substituting 
Eq.~(\ref{comcon}) in ${\cal P}(y,x,l)=\sum_\nu
\frac{y^\nu}{\nu!}{\cal P}_\nu(x,l)$ and 
recognizing the sum as the Taylor expansion 
of the closed form. Binomially expanding 
Eq.~(\ref{confgenf}) yields the polynomial 
formula.~$\bullet$
 
{\bf Corollary~3.4} {\it The partial 
differential equation (PDE) of the 
generating function}
\begin{eqnarray}
\frac{\partial{\cal P}}{\partial x}=-y{\cal 
P}(y,x,l)  
\label{pde}
\end{eqnarray}
{\it is equivalent to the basic differential 
recursion} 
\begin{eqnarray*}
{\cal P'}_\nu(x,l)=-\nu{\cal P}_{\nu-1}(x,l),  
\end{eqnarray*}
{\it and}
\begin{eqnarray}
(1+y)\frac{\partial{\cal P}}{\partial y}=
[c-1+l-x(1+y)]{\cal P}(y,x,l) 
\label{pdee}
\end{eqnarray}
{\it to the three-term recursion} 
\begin{eqnarray*}
(c-1+l-\nu-x){\cal P}_\nu(x,l)=\nu x{\cal 
P}_{\nu-1}(x,l)+{\cal P}_{\nu+1}(x,l). 
\end{eqnarray*}

{\it Proof.} Taking partial derivatives we 
obtain the PDEs~(\ref{pde}), (\ref{pdee}) 
and expanding them in complementary 
polynomials gives the recursions, of which 
the former also follows from Eq.~(32) 
Cor.~4.2~\cite{hjw}.~$\bullet$

{\bf Theorem~3.5} {\it The generating 
function obeys elegant composition laws}
\begin{eqnarray}\nonumber
&&{\cal P}(y,x,l_1){\cal P}(y,x,l_2)={\cal 
P}(y,x,l_1+l_2){\cal P}(y,x,0)\\&&={\cal 
P}(y,x,l_1+l_2-l){\cal P}(y,x,l),~0\le l\le l_1+l_2.
\label{comp}
\end{eqnarray}  
\begin{eqnarray*}
&&\sum_{\nu_1=0}^\nu\left(\nu\atop \nu_1\right)
{\cal P}_{\nu_1}(x,l_1){\cal P}_{\nu-\nu_1}(x,l_2)
\\&&=\sum_{\nu_1=0}^\nu\left(\nu\atop \nu_1\right)
{\cal P}_{\nu_1}(x,l_1+l_2){\cal P}_{\nu-\nu_1}(x,0)
\\&&=\sum_{\nu_1=0}^\nu\left(\nu\atop \nu_1\right)
{\cal P}_{\nu_1}(x,l_1+l_2-l){\cal P}_{\nu-\nu_1}(x,l),~
0\le l\le l_1+l_2.
\end{eqnarray*}
{\it Proof.} The latter follows by expanding 
Eq.~(\ref{comp}) binomially.~$\bullet$ 

{\bf Theorem~3.6}~(Addition Law) 
\begin{eqnarray*}
\sum_{\lambda=0}^\nu\left(c-1+l\atop \nu-\lambda\right)
\frac{\nu!}{\lambda!}{\cal P}_{\lambda}(x_1+x_2,l)=
\sum_{\nu_1=0}^\nu\left(\nu\atop \nu_1\right)
{\cal P}_{\nu_1}(x_1,l){\cal P}_{\nu-\nu_1}(x_2,l)
\end{eqnarray*}
 
\begin{eqnarray*}
(1+y)^{c-1+l}{\cal P}(y,x_1+x_2,l)={\cal P}(y,x_1,l)
{\cal P}(y,x_2,l). 
\end{eqnarray*}
{\it Proof.} The polynomial addition law follows 
similarly from the addition law of the generating 
function.~$\bullet$

The sum on the lhs of the addition law may be 
inverted by expanding similarly 
\begin{eqnarray*}
{\cal P}(y,x_1+x_2,l)=(1+y)^{1-c-l}{\cal P}(y,x_1,l)
{\cal P}(y,x_2,l) 
\end{eqnarray*}
yielding 
\begin{eqnarray*}
{\cal P}_{\nu}(x_1+x_2,l)=\sum_{\nu_j,0\le \nu_1+\nu_2\le \nu}
\left(1-c-l\atop \nu-\nu_1-\nu_2\right)\frac{\nu!}{\nu_1!\nu_2!}
{\cal P}_{\nu_1}(x_1,l){\cal P}_{\nu_2}(x_2,l). 
\end{eqnarray*}

{\bf Theorem~3.7}~{\it The full addition law} 
\begin{eqnarray*}
{\cal P}(y,x_1+x_2,l_1+l_2)=(1+y)^{1-c}
{\cal P}(y,x_1,l_1){\cal P}(y,x_2,l_2)
\end{eqnarray*}
{\it implies similarly}
\begin{eqnarray*}
{\cal P}_{\nu}(x_1+x_2,l_1+l_2)&=&\sum_{\nu_j\ge 
0,0\le \nu_1+\nu_2\le \nu}\left(1-c\atop 
\nu-\nu_1-\nu_2\right)\frac{\nu!}{\nu_1!\nu_2!}\\&\cdot
&{\cal P}_{\nu_1}(x_1,l_1){\cal P}_{\nu_2}(x_2,l_2). 
\end{eqnarray*}

Laguerre and Hermite polynomials are the usual special 
cases of Eq.~(\ref{confl}) including their 
complementary polynomials and their composition 
and addition laws~\cite{hjw},\cite{awh}.  

\section{Hypergeometric ODE} 

For the hypergeometric ODE
\begin{eqnarray}
x(1-x)y''+[c-(a+b+1)x]y'-aby(x)=0,
\label{hyp}
\end{eqnarray}
in the notations of \cite{hjw}, 
\begin{eqnarray*}
&&\sigma(x)=x(1-x),~\tau(x)=c-(a+b+1)x,\\&&
\sigma'=1-2x,\tau-\sigma'=c-1-(a+b-1)x 
\end{eqnarray*}
so that Pearson's ODE for the weight function $w(x)$ 
\begin{eqnarray*}
x(1-x)w'=[c-1-(a+b-1)x]w(x)
\end{eqnarray*}
has the solution 
\begin{eqnarray*}
w(x)=x^{c-1}(1-x)^{a+b-c};~x\in [0,1].
\end{eqnarray*}
In the following we denote the standard regular 
solution of the hypergeometric ODE~(\ref{hyp}) 
simply by $F(a,b;c;x)$ (instead of ${}_2F_1$).   
 
{\bf Definition~4.1} {\it The Rodrigues' formula 
for polynomial solutions $R_l(x)$ of the hypergeometric 
ODE~(\ref{hyp}) is~\cite{hjw} Eq.~(3)}
\begin{eqnarray*}
R_l(x)=x^{1-c}(1-x)^{c-a-b}\frac{d^l}{dx^l}\left(
x^{l+c-1}(1-x)^{l+a+b-c}\right), 
\end{eqnarray*}
{\it with} 
\begin{eqnarray*}
R_l(x)=R_l(0)F(-l,b;c;x),
\end{eqnarray*}
{\it and for the complementary polynomials it 
is~\cite{hjw} Eq.~(5)}
\begin{eqnarray*}
{\cal P}_\nu(x,l)=x^{\nu+1-c-l}(1-x)^{\nu-l+c-a-b}
\frac{d^\nu}{dx^\nu}\left(x^{c-1+l}(1-x)^{l+a+b-c}\right),
\end{eqnarray*}
{\it so that ${\cal P}_l(x,l)=R_l(x).$ The $R_l$ 
also depend on the parameters of the hypergeometric 
ODE, but these are suppressed for simplicity. (The 
$R_l$ here are not to be confused with those in 
Sect.~3.)} 

{\bf Theorem~4.2} {\it The generating function of 
the complementary polynomials is simply a 
polynomial in two variables; it has the closed form} 
\begin{eqnarray}
{\cal P}(y,x,l)=[1+y(1-x)]^{l+c-1}(1-xy)^{l+a+b-c} 
\label{hypgf}
\end{eqnarray}
{\it with} 
\begin{eqnarray*}
{\cal P}_\nu(x,l)=\nu!\sum_{\lambda=0}^\nu
\left(l+a+b-c\atop \lambda\right)\left(l+c-1\atop 
\nu-\lambda\right)(1-x)^{\nu-\lambda}(-x)^\lambda.
\end{eqnarray*}
{\it Proof.} Eq.~(\ref{hypgf}) follows from~\cite{hjw} 
(Theor.3.2). Binomially expanding the generating 
function yields the polynomial expression.~$\bullet$

Note that the generating function~(\ref{hypgf}) is 
just a product of powers of two bilinear forms in 
two variables and amazingly simple, given 
the complexity of the regular solution of 
the hypergeometric ODE. It is certainly much 
simpler than the standard one for Jacobi 
polynomials which, in essence, are the 
complementary polynomials. The latter 
are the main ingredient of the rotation 
matrix elements $d^l_{m m'}(\theta)$ 
commonly used by physicists~\cite{ed}.   
 
Taking partial derivatives leads to 
\begin{eqnarray}\nonumber
&&(1-xy)[1+y(1-x)]\frac{\partial{\cal 
P}}{\partial y}=[(l+c-1)(1-x)(1-xy)\\&&
-x(l+a+b-c)(1+y)(1-x)]{\cal P}(y,x,l)
\label{phyde1}
\end{eqnarray}
and 
\begin{eqnarray}\nonumber 
&&(1-xy)[1+y(1-x)]\frac{\partial{\cal 
P}}{\partial x}=-y[(l+c-1)(1-xy)\\&&
+(l+a+b-c)(1+y)(1-x)]{\cal P}(y,x,l). 
\label{pde2}
\end{eqnarray}
 
{\bf Corollary~4.3} {\it The following 
recursions hold for the complementary 
polynomials} 
\begin{eqnarray}\nonumber
&&{\cal P}_{\nu+1}=[l+c-\nu-1-x(2l-2\nu+a+b-1)]
{\cal P}_\nu\\\nonumber&&-(2l+a+b-1)x(1-x)\nu{\cal 
P}_{\nu-1}+\nu(\nu-1)x(1-x){\cal P}_{\nu-2},\\
\label{rec1}
\end{eqnarray}
\begin{eqnarray}\nonumber
&&(\nu+2l-2){\cal P}_\nu+\nu[(\nu-a-b-2)(1-2x)
+l+a+b-c\\\nonumber&&-x(2l+a+b-1)]{\cal P}_{\nu-1}
-\nu(\nu-1)(\nu-a-b-3)x(1-x)\\&&\cdot{\cal 
P}_{\nu-2}=0. 
\label{rec2}
\end{eqnarray}
{\it Proof.} Expanding the PDE~(\ref{phyde1}) 
in terms of complementary polynomials gives 
the recursion~(\ref{rec1}), while the  
PDE~(\ref{pde2}) yields 
\begin{eqnarray*}
&&{\cal P}'_\nu+\nu(1-2x){\cal P}'_{\nu-1}
-\nu(\nu-1)x(1-x){\cal P}'_{\nu-2}\\&=&
-(2l+a+b-1)\nu{\cal P}_{\nu-1}
-[l+a+b-c-x(2l+a+b-1)]\\&\cdot&\nu(\nu-1)
{\cal P}_{\nu-2}. 
\end{eqnarray*} 
Using the basic recursive ODE~(\ref{badirec}) 
below allows eliminating the derivatives, thus 
yielding the recursion~(\ref{rec2}), which 
reflects the complexity of the regular solution 
of the hypergeometric ODE~(\ref{hyp}).~$\bullet$  

{\bf Theorem~4.4} {\it The following composition 
laws are valid}
\begin{eqnarray*}
&&{\cal P}(y,x,l_1){\cal P}(y,x,l_2)={\cal 
P}(y,x,l_1+l_2){\cal P}(y,x,0)\\&&=
{\cal P}(y,x,l_1+l_2-l){\cal P}(y,x,l),~0\le 
l\le l_1+l_2;  
\end{eqnarray*}  
\begin{eqnarray*}
&&\sum_{\nu_1=0}^\nu\left(\nu\atop \nu_1\right)
{\cal P}_{\nu_1}(x,l_1){\cal P}_{\nu-\nu_1}(x,l_2)
\\&&=\sum_{\nu_1=0}^\nu\left(\nu\atop \nu_1\right)
{\cal P}_{\nu_1}(x,l_1+l_2-l){\cal P}_{\nu-\nu_1}(x,l). 
\end{eqnarray*}
{\it Proof.} Binomially expanding the generating 
function identities yields the polynomial 
composition laws.~$\bullet$

{\bf Theorem~4.5} {\it The ${\cal P}_\nu(x,l)$ 
obey the hypergeometric ODE~(\ref{hyp}) with 
parameters}  
\begin{eqnarray}
(a,b,c)\to (A=-\nu,B=2l-\nu+a+b,C=l-\nu+c)
\label{subst}
\end{eqnarray}  
{\it so that}
\begin{eqnarray*}
{\cal P}_{\nu}(x,l)={\cal P}_{\nu}(0,l)
F(-\nu,B;C;x),~{\cal P}_{\nu}(0,l)=\nu!
\left(l+C-1\atop \nu\right).
\end{eqnarray*}
{\it Proof.} From~\cite{hjw}~(Theor.~5.1) 
the ${\cal P}_{\nu}(x,l)$ obey the ODE  
\begin{eqnarray*}
x(1-x)\frac{d^2}{dx^2}{\cal P}_\nu(x,l)&+&
[(l-\nu)(1-2x)+c-(a+b+1)x]\\&\cdot&\frac{
d{\cal P}_\nu(x,l)}{dx}=-\nu(2l-\nu+a+b)
{\cal P}_\nu(x,l),
\end{eqnarray*}
where the rhs determines $A$ and $B$ in 
Eq.~(\ref{subst}). Then one verifies that 
this ODE is the hypergeometric ODE by 
comparing the coefficient of 
$\frac{d{\cal P}_{\nu}(x,l)}{dx}$ of this 
ODE with Eq.~(\ref{hyp}) after making the 
substitutions (\ref{subst}).~$\bullet$   

The generating function (Theor.~4.2) 
applies to the polynomial solutions of 
Eq.~(\ref{hyp}), albeit in conjunction 
with the parameters of Eq.~(\ref{subst}). 
The complementary polynomials also 
satisfy the basic recursive 
ODE~\cite{hjw} (Cor.~4.2)
\begin{eqnarray}
\frac{d{\cal P}_\nu}{dx}=\nu(\nu-a-b-2){\cal 
P}_{\nu-1}(x,l).
\label{badirec}
\end{eqnarray}

Displaying the parameters, the generating 
function obeys the symmetry 
\begin{eqnarray*}
{\cal P}(-y,1-x,l,c-1,a+b-c)={\cal 
P}(y,x,l,a+b-c,c-1).
\end{eqnarray*}

\section{Relativistic Hermite Polynomials} 

{\bf Definition~5.1} {\it Replacing 
$e^{-x^2}\to (1+\frac{x^2}{N})^{-N},$ 
relativistic Hermite polynomials are 
defined by the Rodrigues' formula~\cite{cv}}
\begin{eqnarray}
H_n^N(x)=(-1)^n\left(1+\frac{x^2}{N}
\right)^{N+n}\frac{d^n}{dx^n}
\left(1+\frac{x^2}{N}\right)^{-N},
\label{dhn}
\end{eqnarray}
which is beyond the framework of ~\cite{hjw}.

{\bf Definition~5.2} {\it The complementary 
polynomials are defined recursively by} 
\begin{eqnarray}
H_n^N(x)=(-1)^{n-\nu}\left(1+\frac{x^2}
{N}\right)^{N+n}\frac{d^{n-\nu}}{dx^{n-\nu}}
\frac{{\cal P}_\nu(x,N)}{\left(1+\frac{x^2}
{N}\right)^{N+\nu}}.
\label{dpn} 
\end{eqnarray}
From $\nu=n$ in Def.~5.2 it follows that 
${\cal P}_\nu(x,N)=H_\nu^N(x),$ just as Hermite 
polynomials coincide with their complementary 
polynomials~\cite{hjw}. Substituting this in 
Eq.~(\ref{dpn}) yields the new relations  
\begin{eqnarray}
H_n^N(x)=(-1)^{n-\nu}\left(1+\frac{x^2}{N}\right)^{N+n}
\frac{d^{n-\nu}}{dx^{n-\nu}}\left(\frac{H_\nu^N(x)}
{(1+\frac{x^2}{N})^{N+\nu}}\right)
\end{eqnarray} 
For $\nu=n-1$ this gives the useful differential 
recursion 
\begin{eqnarray*}
H_n^N(x)-\frac{2x}{N}(N+n-1)H_{n-1}^N(x)+\left(1+
\frac{x^2}{N}\right)\frac{dH_{n-1}^N}{dx}=0. 
\end{eqnarray*}

{\bf Theorem~5.3} {\it The generating function 
$H^N(y,x)$ is given in closed form by} 
\begin{eqnarray*}
H^N(y,x)=\sum_{\nu=0}^\infty\frac{y^\nu}{\nu!}
H_\nu^N(x)=\bigg[\left(1-\frac{xy}{N}\right)^2
+\frac{y^2}{N}\bigg]^{-N}. 
\end{eqnarray*} 
{\it Proof.} This known result follows at once 
and in a new way from using the Rodrigues' 
formula~(\ref{dhn}) for complementary polynomials 
in Theor.~3.2~\cite{hjw}, then recognizing 
the sum as the Taylor expansion of the 
closed form:  
\begin{eqnarray*}
\sum_{\nu=0}^\infty\frac{y^\nu}{\nu!}
H_\nu^N(x)=\left(1+\frac{x^2}{N}\right)^N
\left(1+\frac{X^2}{N}\right)^{-N},~
X=x-y\left(1+\frac{x^2}{N}\right).~\bullet
\end{eqnarray*}

\section{Pre-Laguerre Polynomials}

{\bf Definition~6.1} {\it Replacing the weight function 
$e^{-x}\to (1+\frac{x}{N})^{-N}$ we introduce the 
pre-Laguerre polynomials by the Rodrigues' formula}  
\begin{eqnarray}
{\cal L}_l^N(x)=(1+\frac{x}{N})^{N+l}\frac{d^l}{dx^l}
\frac{x^l}{(1+\frac{x}{N})^N},
\label{rod}
\end{eqnarray}
which is beyond the framework of~\cite{hjw}.

{\bf Proposition~6.2} {\it The complementary 
polynomials ${\cal P}_\nu^N(x,l)$ are defined 
recursively by} 
\begin{eqnarray}
{\cal L}_l^N(x)=(1+\frac{x}{N})^{N+l}\frac{
d^{l-\nu}}{dx^{l-\nu}}\frac{x^{l-\nu}{\cal 
P}_\nu^N(x,l)}{(1+\frac{x}{N})^{N+\nu}}, 
\nu=0,1,\ldots,l.
\label{rodc}
\end{eqnarray}

{\it Proof.} This follows by induction on $\nu$ 
provided they obey the linear differential 
recursion
\begin{eqnarray}\nonumber
{\cal P}_{\nu+1}^N(x,l)&=&[(1+\frac{x}{N})
(l-\nu)-(1+\frac{\nu}{N})x]{\cal P}_\nu^N(x,l)\\&
+&x(1+\frac{x}{N})\frac{d}{dx}{\cal P}_\nu^N(x,l) 
\end{eqnarray}
upon carrying out the innermost differentiation 
in Eq.~(\ref{rodc}).~$\bullet$

For $\nu=0$ we get ${\cal P}_0^N(x,l)=1,$ for 
$\nu=1:~{\cal P}_1^N(x,l)=l(1+x/N)-x;$ for 
$\nu=l:~P_l^N(x)={\cal L}_l^N(x,l).$ 

{\bf Proposition~6.3} {\it The complementary 
polynomials obey their own Rodrigues' formula}  
\begin{eqnarray}
{\cal P}_\nu^N(x,l)=x^{\nu-l}\left(1+\frac{x}{N}
\right)^{N+\nu}\frac{d^\nu}{dx^\nu}\frac{x^l}
{(1+\frac{x}{N})^N}.
\label{rodcp}
\end{eqnarray}

{\it Proof.}  This follows upon plugging 
Eq.~(\ref{rodcp}) into Eq.~(\ref{rodc}) of 
Prop.~6.2, the recursive definition of 
complementary polynomials, implying the 
Rodrigues formula for the ${\cal L}_l^N(x),$ 
Def.~6.1, and vice versa: 
\begin{eqnarray*}
(1+\frac{x}{N})^{N+l}\frac{d^{l-\nu}}{dx^{l-\nu}}
\frac{d^{\nu}}{dx^{\nu}}\frac{x^l}{(1+\frac{x}
{N})^N}={\cal L}_l^N(x).~\bullet
\end{eqnarray*} 
   
More generally we obtain similarly 

{\bf Corollary~6.4} 
\begin{eqnarray*}
{\cal P}_\nu^N(x,l)=x^{\nu-l}(1+\frac{x}{N})^{N+\nu}
\frac{d^{\nu-\mu}}{dx^{\nu-\mu}}\bigg[\frac{x^{l-\mu}
{\cal P}_\mu^N(x,l)}{(1+\frac{x}{N})^{N+\mu}}\bigg]. 
\end{eqnarray*}

{\bf Theorem~6.5} {\it The generating function of 
the complementary polynomials is given in the simple 
closed form by}   
\begin{eqnarray}
{\cal P}^N(y,x,l)=\frac{[1+y(1+\frac{x}{N})]^l}
{(1+\frac{xy}{N})^N}.
\label{genf}
\end{eqnarray}

{\it Proof.} This follows upon applying the 
proof of Theor.~3.2~\cite{hjw} using Prop.~6.3:
\begin{eqnarray*}
&&{\cal P}^N(y,x,l)=\sum_{\nu=0}^\infty\frac{y^\nu}
{\nu!}{\cal P}_\nu^N(x,l)\\&&=\frac{(1+\frac{x}{N})^N}{x^l}
\sum_{\nu=0}^\infty\frac{[xy(1+\frac{x}{N})]^\nu}{\nu!}
\frac{d^\nu}{dx^\nu}\frac{x^l}{(1+\frac{x}{N})^N}\\&&
=(1+\frac{x}{N})^N\frac{[1+y(1+\frac{x}{N})]^l}
{\{1+\frac{x}{N}[1+y(1+\frac{x}{N})]\}^N}, 
\end{eqnarray*}
recognizing the sum as a Taylor expansion, then 
canceling $x^l$ and thereby transforming it into 
Eq.~(\ref{genf}).~$\bullet$

{\bf Proposition~6.6} {\it In general, 
complementary polynomials are given by}
\begin{eqnarray*}
{\cal P}_\nu^N(x,l)=\nu!\sum_{\mu=0}^l
\left(l\atop \mu\right)\left(-N\atop 
\nu-\mu\right)\left(\frac{x}{N}\right)^{
\nu-\mu}(1+\frac{x}{N})^\mu. 
\end{eqnarray*}

{\it Proof.} This follows from expanding 
numerator and denominator powers in 
Eq.~(\ref{rodc}) binomially
\begin{eqnarray*}
\sum_{\nu=0}^\infty\frac{1}{\nu!}
{\cal P}_\nu^N(x,l)=\bigg[\sum_{\mu=0}^l
\left(l\atop \mu\right)y^\mu(1+\frac{x}
{N})^\mu\bigg]\bigg[\sum_{\lambda=0
}^\infty\left(-N\atop \lambda\right)
\left(\frac{xy}{N}\right)^\lambda\bigg]
\end{eqnarray*} 
and collecting terms with $\nu=\lambda+\mu.~\bullet$

{\bf Proposition~6.7} {\it A recursive formula 
for the complementary polynomials is} 
\begin{eqnarray}
\sum_{\lambda=0}^\nu\left(N\atop \nu-\lambda\right)
\left(\frac{x}{N}\right)^{\nu-\lambda}
\frac{{\cal P}_\lambda^N(x,l)}{\lambda!}=\cases{
\left(l\atop \nu\right)(1+\frac{x}{N})^\nu,\nu\le l;\cr
&\cr
0,\nu>l.\cr}
\end{eqnarray}

{\it Proof.} This follows by expanding binomially 
the numerator and denominator of Eq.~(\ref{genf}) 
in the form 
\begin{eqnarray*}
(1+\frac{xy}{N})^N{\cal P}^N(y,x,l)=
[1+y(1+\frac{x}{N})]^l
\end{eqnarray*}
of Eq.~(\ref{genf}).~$\bullet$ 

Taking derivatives and then expanding binomially 
the form 
\begin{eqnarray*}
(1+\frac{xy}{N})^N\frac{\partial{\cal P}^N(y,x,l)}
{\partial x}&=&\frac{ly}{N}(1+\frac{x}{N})
[1+y(1+\frac{x}{N})]^{l-1}\\&-&y[1+y(1+\frac{x}{N})]^l
[1+y(1+\frac{x}{N})]^l
\end{eqnarray*}
of Eq.~(\ref{genf}) yields a general formula 
and recursions for the derivatives of the 
complementary polynomials. 

Upon applying Pearson's ODE 
\begin{eqnarray*}
(\sigma w)'=\tau w,~\sigma=x,~w=(1+\frac{x}{N})^{-N}
\end{eqnarray*}
gives $\tau(x)=1-\frac{x}{1+\frac{x}{N}}.$ For 
$N\to\infty$ the standard Laguerre ODE form 
results, $\tau=1-x.$ Thus, ODEs for 
${\cal L}_l, {\cal P}_\nu$ will take the forms 
\begin{eqnarray*}
x\frac{d^2}{dx^2}{\cal L}_l^N+\left(1-\frac{x}
{1+\frac{x}{N}}\right)\frac{d{\cal L}_l^N}{dx}
+V_l^N(x){\cal L}_l^N
=\lambda_l {\cal L}_l^N
\end{eqnarray*} 
\begin{eqnarray*}
x\frac{d^2}{dx^2}{\cal P}_\nu^N+\left(1-\frac{x}
{1+\frac{x}{N}}\right)\frac{d{\cal P}_\nu^N}{dx}
+v_\nu^N(x){\cal P}_\nu^N=\Lambda_\nu {\cal P}_\nu^N. 
\end{eqnarray*} 
If such ODEs are to be the radial Schr\"odinger 
equation of a quantum mechanical problem, then 
the derivative term $(x-\frac{x}{1+\frac{x}{N}})
\frac{d}{dx}$ has to be considered as part of 
the potential making it non-local. 

{\bf Corollary~6.8} (Recursion) {\it The PDE} 
\begin{eqnarray}\nonumber
&&\left(1+\frac{xy}{N}\right)[1+y(1+\frac{x}{N})]
\frac{\partial{\cal P}^N}{\partial y}\\&&=
\{l(1+\frac{x}{N})(1+\frac{xy}{N}-x[1+y(1+
\frac{x}{N})]\}{\cal P}^N(y,x,l)
\label{pde1}
\end{eqnarray}
{\it is equivalent to the recursion}
\begin{eqnarray}\nonumber
&&{\cal P}^N_{\nu+1}+[\nu-l+x(\frac{2\nu-l}{N}+1)]
{\cal P}^N_\nu\\&&+\nu(\frac{\nu-l-1}{N}+1)x(1+
\frac{x}{N}){\cal P}^N_{\nu-1}=0.
\label{rec3} 
\end{eqnarray}  
{\it Proof.} Taking the partial derivative of 
${\cal P}^N$ yields the PDE~(\ref{pde1}),  
and expanding it in complementary polynomials 
gives the recursion~(\ref{rec3}).~$\bullet$
 
The elegant angular momentum addition identity 
\begin{eqnarray*}
{\cal P}^N(y,x,l_1+l_2)=(1+\frac{xy}{N})^N
{\cal P}^N(y,x,l_1){\cal P}^N(y,x,l_2)
\end{eqnarray*}
and in the form 
\begin{eqnarray*}
{\cal P}^N(y,x,l_1){\cal P}^N(y,x,l_2)=
(1+\frac{xy}{N})^{-N}{\cal P}^N(y,x,l_1+l_2)
\end{eqnarray*}
yield, upon expanding the generating functions 
into their complementary polynomials: 

{\bf Theorem~6.9} (Composition Laws)  
\begin{eqnarray*}
&&{\cal P}_\nu^N(x,l_1+l_2)=\sum_{\mu=0}^\nu\frac{
\nu!}{\mu!}\left(N\atop \nu-\mu\right)\left(
\frac{x}{N}\right)^{\nu-\mu}\\&&\cdot\sum_{\nu_1=0}^\mu
\left(\mu\atop \nu_1\right){\cal P}_{\nu_1}^N(x,l_1)
{\cal P}_{\mu-\nu_1}^N(x,l_2) 
\end{eqnarray*} 
{\it and}
\begin{eqnarray*}\nonumber
&&\sum_{\nu_1=0}^\nu\left(\nu\atop \nu_1\right)
{\cal P}_{\nu_1}^N(x,l_1){\cal P}_{\nu-\nu_1}^N(x,l_2)
\\&&=\sum_{\lambda=0}^\nu\frac{\nu!}{\lambda!}\left(
-N\atop \nu-\lambda\right)\left(\frac{x}{N}\right)^{
\nu-\lambda}{\cal P}_\lambda^N(x,l_1+l_2).
\end{eqnarray*}

Finally we get to the special values. For $x=0$ 
in Eq.~(\ref{genf}) ${\cal P}^N(y,0,l)=(1+y)^l$ 
of Theor.~6.4 we find 
\begin{eqnarray*}
{\cal P}_\nu^N(0,l)=\cases{\nu!\left(l\atop 
\nu\right),~\nu\le l;\cr
&\cr
0,\nu>l;\cr}
\end{eqnarray*}  
while for $x=-N$ we find from ${\cal P}^N(y,-N,l)
=(1-y)^{-N}$
\begin{eqnarray*}
{\cal P}_\nu^N(-N,l)=\nu!(-1)^\nu\left(-N\atop 
\nu\right). 
\end{eqnarray*}

\section{More General Framework}
 
Summarizing and generalizing Sects.~5 and 6, 
for given polynomials $\sigma(x), w(x)$ the 
Rodrigues' relation 
\begin{eqnarray}
S_l(x)=w(x)^{l+N}\frac{d^l}{dx^l}\frac{
\sigma(x)^l}{w(x)^N}
\label{rodb}
\end{eqnarray}
defines polynomials $S_l$ with $S_0=1.$ 
Complementary polynomials are defined 
recursively by
\begin{eqnarray}
S_l(x)=w(x)^{l+N}\frac{d^{l-\nu}}{dx^{l-\nu}}
\frac{\sigma^{l-\nu}{\cal P}_\nu(x,l)}{w^{N+\nu}},
\label{cprd}
\end{eqnarray}
which satisfy their Rodrigues' relation
\begin{eqnarray}
{\cal P}_\nu(x,l)=\sigma^{\nu-l}w^{N+\nu}
\frac{d^\nu}{dx^\nu}\frac{\sigma^l}{w^N}. 
\label{rocp}
\end{eqnarray}
Substituting Eq.~(\ref{rocp}) in Eq.~{\ref{cprd}} 
reproduces Eq.~(\ref{rodb}). Using Eq.~(\ref{rocp}) 
in the generating function sum yields its closed 
form
\begin{eqnarray*}
{\cal P}(y,x,l)=\left(\frac{w(x)}{w(x+y\sigma(x)
w(x))}\right)^N\left(\frac{\sigma(x+y\sigma(x)
w(x))}{\sigma(x)}\right)^l. 
\end{eqnarray*}  

Generalizing similarly Sect.~4 on the 
hypergeometric ODE, the Rodriguez' 
formula 
\begin{eqnarray}
S_l(x)=\frac{1}{w^a(x)\sigma^b(x)}\frac{d^l}
{dx^l}(w^{a+l}\sigma^{b+l})
\label{gyp}
\end{eqnarray}
defines polynomials $S_l$ and their 
complementary set
\begin{eqnarray*}
{\cal P}_\nu(x,l)=\frac{1}{w^{a+l}\sigma^{b+l}}
\frac{d^{\nu}}{dx^{\nu}}(w^a\sigma^b).
\end{eqnarray*}

For $w(x)=x,\sigma(x)=1-x$ and appropriate 
parameters $a,b$ these formulas reduce to 
those of the hypergeometric ODE. The 
$l$-dependence $(w\sigma)^l$ in Eq.~(\ref{gyp}) 
can be removed by differentiating 
$w^a\sigma^b S_l~\lambda$-fold  
\begin{eqnarray*}
S_{l+\lambda}(x)=\frac{1}{w^a\sigma^b}\frac{
d^\lambda}{dx^\lambda}(w^a\sigma^b S_l(x),
\end{eqnarray*} 
yielding an index translation formula, which 
may be resummed as a translation formula
\begin{eqnarray*}
&&w^a(x+h)\sigma^b(x+h)S_l(x+h)=\\&&\sum_{
\lambda=0}^\infty\frac{h^\lambda}{\lambda!}
\frac{d^\lambda}{dx^\lambda}(w^a\sigma^b S_l(x))
=w^a\sigma^b\sum_{\lambda=0}^\infty S_{l+\lambda}(x).
\end{eqnarray*} 

The generating function ${\cal S}(y,x,l)$ for 
the complementary polynomials of the $S_l$ is 
\begin{eqnarray*}
{\cal S}(y,x,l)=\left(\frac{w(x+y\sigma(x)w(x))}
{w(x)}\right)^{a+l}\left(\frac{\sigma(x+y\sigma(x)
w(x))}{\sigma(x)}\right)^{b+l}.
\end{eqnarray*}
The complementary polynomials are defined as
\begin{eqnarray*}
{\cal S}_\nu(x,l)=\frac{1}{w^{a-\nu}\sigma^{b-\nu}}
\frac{d^\nu}{dx^\nu}(w^{a+l}\sigma^{b+l}).
\end{eqnarray*}
For Legendre polynomials $w(x)=1+x,\sigma(x)=1-x,
a=b=0$ and
\begin{eqnarray*}
P_l(x+h)=\sum_{\lambda=0}^\infty(-2h)^\lambda
\left(l+\lambda\over l\right)P_{l+\lambda}(x). 
\end{eqnarray*}

Starting from a pair of polynomials this 
generalized framework provides sets 
polynomials and their complementary sets  
that are summed up in closed form.  

\section{Summary And Conclusions} 

The recursion relation for the normalizations 
of the complementary Legendre polynomials 
completes their theory~\cite{trio}. 

For the confluent and hypergeometric ODEs 
the result of central importance is that 
their complementary polynomials obey the 
same ODE, albeit with more general 
parameters, so that their generating 
functions in their simple closed form 
contain them all.



\begin{thebibliography}{0}  

\bibitem{hjw} H. J. Weber, {\it Connections                                     
between real polynomial solutions of                                            
hypergeometric-type differential equations                                      
with Rodrigues formula,} Central Eur. J.
Math. {\bf 5}(2)(2007) 415-427.

\bibitem{awh} G. B. Arfken, H. J. Weber, 
F. Harris, {\it Mathematical Methods for 
Physicists}, 7th ed., Elsevier, Amsterdam 
(2013). 

\bibitem{trio} H. J. Weber, {\it Sonata For 
Trio: Associated Legendre Functions from 
Rodrigues' Formula for Legendre Polynomials}, 
Int. J. Appl. Math. Mech. 3(1) (2014) 9-20. 

\bibitem{rom} H. J. Weber, {\it Connections 
between Romanovski and other polynomials}, 
Central Eur. J. Math. 5(3)(2007) 581-595; 

\bibitem{rwak} A. P. Raposo, H. J. Weber, 
D. E. Alvarez--Castillo, M. Kirchbach, Central 
Eur. J. Phys. 5(3) (2007) 253-284; https://en. 
wikipedia.org/wiki/Romanovskipolynomials.  

\bibitem{ed} A. R. Edmonds, {\it Angular Momentum 
In Quantum Mechanics}, Princeton Univ. Press (1957). 

\bibitem{cv} C. Vignat, {\it More Old And New 
Results About Relativistic\\Hermite Polynomials}, 
J. Math. Phys. 52(9)(2011) 093503, 16 p. and 
Refs. therein.
 
\end{thebibliography}
\end{document}